\documentclass[10pt]{amsart}
\usepackage{mathptmx}
\usepackage{amsmath}
\usepackage{amssymb}
\usepackage{array}
\usepackage{geometry}
\usepackage[bookmarks=true,colorlinks=true, pdfstartview=FitV, linkcolor=black, citecolor=blue, urlcolor=black]{hyperref}

\usepackage{color}
\definecolor{DarkRed}{rgb}{0.55,.00,0.2}
\definecolor{DarkGrey}{rgb}{0.35,.35,0.35}

\theoremstyle{definition}

\theoremstyle{remark}

\numberwithin{equation}{section}

\begin{document}

\title{ Towards  the Casas-Alvero conjecture}

\author{S. Yakubovich}
\address{Department of Mathematics, Fac. Sciences of University of Porto,Rua do Campo Alegre,  687; 4169-007 Porto (Portugal)}
\email{ syakubov@fc.up.pt}


\date{\today}


\maketitle

\markboth{\rm \centerline{ S.  YAKUBOVICH}}{}
\markright{\rm \centerline{CASAS-ALVERO CONJECTURE}}

\begin{abstract}  We investigate necessary and sufficient conditions for an arbitrary polynomial of degree $n$ to be trivial, i.e. to have the form $a(z-b)^n$.  These results are related to an open   problem, conjectured in 2001  by   E. Casas- Alvero.  It says,    that  any complex univariate polynomial, having a common root with each of its non-constant derivative  must be a power of a linear polynomial.  In particular, we establish  determinantal representation  of the Abel-Goncharov interpolation polynomials,  related to the problem and having its  own interest.   Among other results are new  Sz.-Nagy type identities  for complex roots and a generalization of the Schoenberg conjectured analog of Rolle's theorem for polynomials with real and complex coefficients. 

\end{abstract}

\maketitle

\vspace{1cm}

In 2001  E. Casas-Alvero   \cite{CA} conjectured that an arbitrary polynomial $f$  degree $n
\ge 1$ with complex coefficients  of  degree $n \in \mathbb{N}$
$$f(z)= a_0 z^{n}+ a_1 z^{n-1}+ \dots + a_{n-1} z + a_n, \ a_0\neq 0 \eqno(1)$$
is of the trivial monomial form   $f(z)= a(z-b)^n, a, b \in \mathbb{C}$,    if   and \ only   if $f$  shares a root with each of its derivatives $f^{(1)}, f^{(2)}, \dots,   f^{(n-1)}.$  It is proved for small degrees, for infinitely many degrees, for instance, for all powers $n$,   when $n$ is a prime   (see in \cite{Drai}, \cite{Graf} ).  We will call these common roots of the corresponding derivatives by $z_1,  z_{2}, \dots, z_{n-1} \in \mathbb{C}$ (repeated terms are permitted).  As it was recently observed by the author \cite{jca}, the polynomial (1), satisfying the Casas- Alvero conditions (the CA-polynomial),  can be identified, involving  the familiar Abel-Goncharov interpolation polynomials \cite{eva},  which are defined by the following recurrence relation
$$
     f(z)\equiv  G_{n}(z, z_{0}, z_{1},\dots, z_{n-1})  = z^n  -  \sum_{k=0}^{n-1} {{n}\choose {k}}   z_k ^{n- k} 
     G_{k}(z, z_{0}, z_{1},\dots, z_{k-1}),\  G_0(z)\equiv 1\eqno(2)
$$
with the additional conditions ($z_{0}$ is a root of $f$)
$$G_{n}(z_j , z_{0}, z_{1},\dots, z_{n-1}) =0,\quad j= 1,2,\dots, n-1. \eqno(3)$$
It is known \cite{eva}, that the Abel-Goncharov polynomial can be represented as a multiple integral in the complex plane 
$$f(z)= G_{n}(z, z_{0}, z_{1},\dots, z_{n-1}) = n! \int_{z_{0}}^{z} \int_{z_{1}}^{s_{1}} \dots \int_{z_{n-1}}^{s_{n-1}} ds_{n}\dots ds_{1}\eqno(4)$$
Moreover, making  simple changes  of variables in (4),  it can be verified that $G_n(z)$ is shift-invariant and 
 a homogeneous function of degree $n$ (cf. \cite{leva}).  Namely, for any $\alpha \in \mathbb{C} \backslash \{0\},\ \beta \in \mathbb{C}$   it has 
$$G_n(\alpha z+\beta ) \equiv  G_n\left(\alpha z+\beta, \alpha z_0+\beta, \alpha z_1+\beta,\dots, \alpha z_{n-1}+\beta\right)  = \alpha^n G_{n}(z, z_{0}, z_{1},\dots, z_{n-1}) \equiv  \alpha^n G_n(z).\eqno(5)$$
Without loss of generality one can assume in the sequel that $f$ is a monic polynomial of degree $n\in \mathbb{N} $, i.e. $a_0=1$ in (1).  Generally, it has $k$ distinct roots $\lambda_j$  of multiplicities $r_j, \   j= 1,  \dots,  k, 1\le k\le n$ such that
$$ r_1+ r_2+ \dots r_{k} = n.\eqno(6)$$
By $r$ we denote the maximum of multiplicities in (6), i.e. $r=\max_{1\le j\le k} r_j$.  Since for $n=1$ the polynomial is trivial, we will consider $n \ge 2$.  Moreover, a possible non-trivial CA-polynomial cannot have all distinct roots,  because at least one root is common with its first derivative. Therefore the maximum of multiplicities is at least $2$ and a maximum of possible distinct roots is $n-1$.   Another observation tells that a polynomial whose distinct roots are of the same multiplicity $m\ge 2$, i.e.
$$f(z)= \left[ (z-\lambda_1) (z-\lambda_2)\dots (z-\lambda_k)\right] ^m$$
of degree $n=km$ cannot satisfy the Casas-Alvero conditions since the derivative $f^{(m)} $ has no common roots with $f$.  Consequently, at least two roots are of  different multiplicities. 

Recently (see \cite{jca}), the author proved the following propositions.

{\bf Proposition 1.}  {\it  A polynomial with only real roots of degree $n\ge 2$  is trivial, if and only if its $n-2$nd derivative has a double root.}

{\bf Proposition 2.}  {\it  A possible non-trivial polynomial $f$ of degree $n \ge  6$ with only real roots, sharing a root with  its $n-2$nd and $n-1$st derivatives,  has at least five distinct roots}.

{\bf Proposition 3.}  {\it  A possible non-trivial polynomial $f$ of degree $n \ge 7$ with only real roots, sharing  roots with  its $n-2$nd and $n-1$st derivatives, where roots of the $n-2$nd derivative have different multiplicities as roots of $f$,  has at least six  distinct roots}.

Basing on the homogeneity property (5) we also proved

{\bf Proposition 4.}  {\it  The Casas-Alvero conjecture holds true, if and only if  it is true for common roots lying in the unit circle. }

{\bf Proposition 5.}  {\it  A possible non-trivial CA-polynomial  with only real zeros has at least $5$ distinct roots. }

Concerning the Abel-Goncharov polynomials  it was shown in \cite{jca} that $G_{n}$ satisfies the following upper bound 
$$\left| G_{n}(z, z_{0}, z_{1},\dots, z_{n-1})\right| \le \sum_{k_{0}=0}^{1} \sum_{k_{1}=0}^{2-k_{0}}
\dots \sum_{k_{n-2}=0}^{n-1-k_{0}-k_{1}- \dots -k_{n-3}} { {n}\choose {k_{0}, k_{1},\dots,  k_{n-2},
\  n- k_{0} -k_{1}- \dots- k_{n-2}} }$$ 
$$\times \prod_{s=0}^{n-1} \left| z_{n-2-s} - z_{n-1-s}\right|^{k_{s}},$$
where $z_{-1}\equiv z$ and 
$${ {n}\choose {l_{0}, l_{1},\dots,  l_{m} }} = \frac{n!} { l_{0} ! l_{1}! \dots l_{m} !},\quad  l_{0}+ l_{1}\dots + l_{m} =n.$$
This estimate is sharper than the classical Goncharov upper bound \cite{eva}
$$\left| G_{n}(z, z_{0}, z_{1},\dots, z_{n-1})\right| \le \left( \left|z-z_{0}\right| + \sum_{s=0}^{n-2}
|z_{s+1} - z_{s}| \right)^{n}.$$
These  polynomials can be represented via the so-called Levinson binomial type expansion (see in \cite{leva}, p. 732) 
$$ G_{n}(z, z_{0}, z_{1},\dots, z_{n-1}) = \sum_{k=1}^n (z^k-z_0^k ) {{n}\choose {k}} H_{n-k},\eqno(7)$$
where $H_0=1$ and 
$$H_{n-k} \equiv H_{n-k} (z_k, z_{k+1},\dots, z_{n-1}) = (n-k)!  \int_{z_{k}}^{0} \int_{z_{k+1}}^{s_{k+1}} \dots \int_{z_{n-1}}^{s_{n-1}} ds_{n}\dots ds_{k+1},\  k= 1,\dots, n-1.\eqno(8)$$
However we  will prove in turn that  the polynomials $H_{n-k}$  can be represented in a  determinantal  form of an upper Hessenberg matrix $(n-k)\times (n-k)$ with the entries equal to $1$ on the main subdiagonal \cite{horn}. 

Generally, it has 

{ \bf Lemma 1.}  {Let $n \in \mathbb{N}$ and $a_j \in \mathbb{C},\ j=  1,\dots, n$.  Then }
$${(-1)^{n} \over n!}   H_{n} (a_1,  a_{2},\dots,  a_{n})  =  
\begin{vmatrix}

a_1& {a_1^2\over 2!} & {a_1^3\over 3!} & \dots&    {a_1^{n} \over n!} \\

1&  a_{2} & {a_{2}^2\over 2!}&  \dots&    {a_{2}^{n-1} \over (n-1)!} \\

0& 1&  a_{3} &   \dots&    {a_{3}^{n-2} \over (n-2)!} \\
 
  \vdots &  \dots& \ddots&    \ddots&    \vdots \\

0& 0& \dots& 1&       a_{n} \\

\end{vmatrix}
 .\eqno(9)$$

\begin{proof}  Appealing to  the principle of mathematical induction and  easily verifying  formula (9) for $n=1, 2$ via the calculation of  the corresponding integral (8),  i.e.  
$$H_1(a_1)= - a_1,\quad   {1 \over 2} \   H_{2} (a_1,  a_{2})  =  
\begin{vmatrix}

a_1& {a_1^2\over 2!}  \\

1&  a_{2}  \\
\end{vmatrix}
=  a_1 a_2-  {a_1^2\over 2!}, $$
we assume that the statement holds for all $1\le k \le n$ and will prove it for $n+1$.  Indeed, expanding the corresponding determinant along the first column by the Laplace theorem,  we find 
$$\begin{vmatrix}

a_1& {a_1^2\over 2!} & {a_1^3\over 3!} & \dots&    {a_1^{n+1} \over (n+1)!} \\

1&  a_{2} & {a_{2}^2\over 2!}&  \dots&    {a_{2}^{n} \over n!} \\

0& 1&  a_{3} &   \dots&    {a_{3}^{n-1} \over (n-1)!} \\
 
  \vdots &  \dots& \ddots&    \ddots&    \vdots \\

0& 0& \dots& 1&       a_{n+1} \\

\end{vmatrix}=  a_1 \begin{vmatrix}

 a_{2} & {a_{2}^2\over 2!}&  \dots&    {a_{2}^{n} \over n!} \\

1&  a_{3} &   \dots&    {a_{3}^{n-1} \over (n-1)!} \\
 
  \dots& \ddots&    \ddots&    \vdots \\

 0& \dots& 1&       a_{n+1} \\

\end{vmatrix} \   - \   \begin{vmatrix}

 {a_1^2\over 2!} & {a_1^3\over 3!} & \dots&    {a_1^{n+1} \over (n+1)!} \\

1&  a_{3} &   \dots&    {a_{3}^{n-1} \over (n-1)!} \\
 
   \dots& \ddots&    \ddots&    \vdots \\

 0& \dots& 1&       a_{n+1} \\

\end{vmatrix}. \eqno(10)$$
The first $n\times n$ determinant in the right-hand side of the latter equality by the induction hypothesis is equal to 
$$\begin{vmatrix}

 a_{2} & {a_{2}^2\over 2!}&  \dots&    {a_{2}^{n} \over n!} \\

1&  a_{3} &   \dots&    {a_{3}^{n-1} \over (n-1)!} \\
 
  \dots& \ddots&    \ddots&    \vdots \\

 0& \dots& 1&       a_{n+1} \\

\end{vmatrix} =  {(-1)^{n} \over n!}   H_{n} (a_2,  a_{3},\dots,  a_{n+1}) .$$
The second determinant we will expand in the same fashion  to obtain
$$ \begin{vmatrix}

 {a_1^2\over 2!} & {a_1^3\over 3!} & \dots&    {a_1^{n+1} \over (n+1)!} \\

1&  a_{3} &   \dots&    {a_{3}^{n-1} \over (n-1)!} \\
 
   \dots& \ddots&    \ddots&    \vdots \\

 0& \dots& 1&       a_{n+1} \\

\end{vmatrix} =   {a_1^2\over 2!}\  \begin{vmatrix}

  a_{3} &    {a_{3}^{2} \over 2!}&   \dots&    {a_{3}^{n-1} \over (n-1)!} \\
 
 1&    a_4&    \dots&  {a_{4}^{n-2} \over (n-2)!}\\

\dots& \ddots&    \ddots&    \vdots \\

0&  \dots& 1&       a_{n+1} \\

\end{vmatrix}-  \begin{vmatrix}

 {a_{1}^{3} \over 3!}&    {a_{1}^{4} \over 4!}& \dots&    {a_{1}^{n+1} \over (n+1)!} \\
 
 1&    a_4&    \dots&  {a_{4}^{n-2} \over (n-2)!}\\

\dots& \ddots&    \ddots&    \vdots \\

0&  \dots& 1&       a_{n+1} \\

\end{vmatrix}. $$
Continuing the same process and applying every time the induction hypothesis, we arrive at the final expansion of the determinant in the left-hand side of the equality (9) 
$$ (-1)^{n+1}  \begin{vmatrix}

a_1& {a_1^2\over 2!} & {a_1^3\over 3!} & \dots&    {a_1^{n+1} \over (n+1)!} \\

1&  a_{2} & {a_{2}^2\over 2!}&  \dots&    {a_{2}^{n} \over n!} \\

0& 1&  a_{3} &   \dots&    {a_{3}^{n-1} \over (n-1)!} \\
 
  \vdots &  \dots& \ddots&    \ddots&    \vdots \\

0& 0& \dots& 1&       a_{n+1} \\

\end{vmatrix} = -   \sum_{k=1} ^{n+1} {a_1^k\over k!} \   H_{n+1-k } (a_{k},  a_{3},\dots,  a_{n+1}).\eqno(11) $$
Fortunately,  the expression in the right-hand side of (8) is calculated by Levinson via the Taylor theorem ( cf.  \cite{leva}, p. 731)  and we find 
$$ -   \sum_{k=1} ^{n+1} {a_1^k\over k!} \   H_{n+1-k } (a_{k},  a_{3},\dots,  a_{n+1}) =  {1\over (n+1)!} 
\ H_{n+1 } (a_{1},  a_{2},\dots,  a_{n+1}).$$
Thus we get the validity of equality (9) for all $n \in \mathbb{N}$ and complete the proof of the Lemma.  
\end{proof}

{\bf Corollary 1.} {\it The Levinson polynomials $H_n,\  n \in \mathbb{N}$  have the following determinantal  representation}

$$  H_{n} (a_1,  a_{2},\dots,  a_{n})  =   (-1)^{n}  
\begin{vmatrix}

a_1&  \  a_1^2&   \  a_1^3&  \dots&    \  a_1^{n} \\

1&    {{2}\choose {1}}  a_{2} &    {{3}\choose {1}}  a_{2}^2&  \dots&    {{n}\choose {1}}   a_{2}^{n-1}  \\

0& 1&   {{3}\choose {2}}  a_{3} &   \dots&    {{n}\choose {2}}   a_{3}^{n-2}  \\

  \vdots &  \dots& \ddots&    \ddots&    \vdots \\

0& 0& \dots& 1&      {{n}\choose {n-1}}   a_{n} \\

\end{vmatrix}
 .\eqno(12)$$
 
 \begin{proof} In fact, the proof  easily follows from (9). For this  we multiply the last column of the determinant by $n!$, the second  column by $2!$, the third one by $3!$  etc.,  and the $n-1$-th column by $(n-1)!$. Then,  dividing the third  row  by $2!$, the fourth row by $3!$ etc., and the last row by $(n-1)!$, we get the result. 
 
  \end{proof} 
  
  The determinantal form (12) can be involved to investigate the Casas-Alvero conjecture.  Precisely, the shift-invariant property (5) for the Abel-Goncharov polynomials allows to suppose without loss of generality   that one of the polynomial roots, say, $z_0=0$.   Then for a fixed sequence $\{z_j\}_1^{n-1}$ of common roots of $f$ and its derivatives up to the  order $n-1$ (7) implies 
$$ G_{n}(z, 0, z_{1},\dots, z_{n-1}) = \sum_{k=1}^n z^k  {{n}\choose {k}} H_{n-k}.\eqno(13)$$
Assuming the existence of a possible non-trivial CA-polynomial $f$, it follows that  at least one of  $z_{1},\dots, z_{n-1}$ is nonzero,  otherwise the polynomial has the   unique root of  multiplicity $n$, which is equal to zero.   Let  $f$   have $s$ nonzero common roots $(1 \le s \le n-1)$ in our sequence  $\{z_j\}_1^{n-1}$
$$  z_j,\  j = i_1, \dots,  i_s,\    \{ i_1, \dots,  i_s \}  \subset  \{1,\dots, n-1\},\   1\le i_1< i_2< \dots < i_s \le n-1,  \eqno(14)$$
corresponding  the derivatives of the order $i_1,\dots, i_s$, respectively.   Hence, appealing to  Corollary 1 and observing from (4) and (8) that the left-hand side of (13) equals to $- H_{n}(z,  z_{1},\dots, z_{n-1})$, we end up with the equality 
$$  H_{n} (z, z_1,  z_{2},\dots,  z_{n-1})  =   (-1)^{n}  
\begin{vmatrix}

z &  \  z^2&   \  z^3&  \dots&    \  z^{n} \\

1&    {{2}\choose {1}}  z_{1} &    {{3}\choose {1}}  z_{1}^2&  \dots&    {{n}\choose {1}}   z_{1}^{n-1}  \\

0& 1&   {{3}\choose {2}}  z_{2} &   \dots&    {{n}\choose {2}}   z_{2}^{n-2}  \\

  \vdots &  \dots& \ddots&    \ddots&    \vdots \\

0& 0& \dots& 1&      {{n}\choose {n-1}}   z_{n-1} \\

\end{vmatrix}
 .\eqno(14)$$
Moreover,  the determinant (14) can be compressed, eliminating   rows,   containing only one nonzero element of the main subdiagonal of Hessenberg's matrix, which is equal to $1$.  Therefore,   we arrive at the equality for any $z\neq 0$
$$  (-1)^{s} \  z^{- i_1} H_{n} (z, z_1,  z_{2},\dots,  z_{n-1})  =   
\begin{vmatrix}

1 &     z^{i_2-i_1}&  \dots&    z^{i_s-i_1}&     \  z^{n-i_1} \\

1&    {{i_2}\choose {i_1}}  z^{i_2-i_1} _{i_1} &   \dots&   {{i_k}\choose {i_1}}  z_{i_1}^{i_s-i_1}&    {{n}\choose {i_1}}   z_{i_1}^{n-i_1}  \\

0& 1&   {{i_3}\choose {i_2}}  z_{i_2}^{i_3- i_2} &  \dots&     {{n}\choose {i_2}}   z_{i_2}^{n-i_2}  \\

  \vdots &  \dots& \ddots&    \ddots&    \vdots \\

0& 0& \dots& 1&      {{n}\choose {i_s}}   z_{i_s}^{n-i_s} \\

\end{vmatrix}
 ,\eqno(15)$$
containing the determinant of the order $(s+1)\times (s+1)$.   Now, taking into account the conditions (see (3))
$$H_{n}(z_j,  z_{1},\dots, z_{n-1})  =0, \  z_j\neq 0, \quad j= i_1, i_2,\dots, i_s,\eqno(16)$$
we have 
$$  
\begin{vmatrix}

1 &     z_j^{i_2-i_1}&  \dots&    z_j^{i_s-i_1}&     \  z_j^{n-i_1} \\

1&    {{i_2}\choose {i_1}}  z^{i_2-i_1} _{i_1} &   \dots&   {{i_k}\choose {i_1}}  z_{i_1}^{i_s-i_1}&    {{n}\choose {i_1}}   z_{i_1}^{n-i_1}  \\

0& 1&   {{i_3}\choose {i_2}}  z_{i_2}^{i_3- i_2} &  \dots&     {{n}\choose {i_2}}   z_{i_2}^{n-i_2}  \\

  \vdots &  \dots& \ddots&    \ddots&    \vdots \\

0& 0& \dots& 1&      {{n}\choose {i_s}}   z_{i_s}^{n-i_s} \\

\end{vmatrix} =0,\quad  j=  i_1,\dots, i_s.\eqno(17)$$
Let $s=1$, i.e. the polynomial has only one nonzero  root $z_{i_1} \neq 0$ in our sequence of common roots $z_1,z_2, \dots, z_{n-1}$, which corresponds to the $i_1$-th derivative.  Hence,  putting in (17) $j=i_1$,  the determinant 
$$  
\begin{vmatrix}

1 &       \  z_{i_1} ^{n-i_1} \\
1&      {{n}\choose {i_1}}   z_{i_1} ^{n-i_1} \\
\end{vmatrix} =0.\eqno(18)$$
But this is impossible, since  we have $ {{n}\choose {i_1}}  =1,\ n > i_1\ge1.$    Hence $2 \le s \le n-1$.   Returning to (13) and using conditions (16), we sum up these $s$ equalities to obtain 

$$ \sum_{j=1}^{s+1}    \sum_{m=1}^{s+1} z_{i_j} ^{i_m}  {{n}\choose {i_m}} H_{n-i_m}(z_{i_m}, 0,\dots, z_{i_{m+1}},\dots, z_{n-1})  =  \sum_{m=1}^{s+1}   {{n}\choose {i_m}}  P_m H_{n-i_m} = 0,\eqno(19)  $$
where $i_{s+1}= n$ and $P_m=  \sum_{j=1}^{s+1}   z_{i_j}^{i_m}$.    Equalities (19) are analogs of the familiar Newton identities  for nonzero roots of a possible non-trivial CA-polynomial.

Recalling distinct roots of $f(z)$ $\lambda_j,\ j=1,\dots, k$ and assuming that the root of the $n-1$st derivative takes, for instance, the value $\lambda_1$, we use the first Vi\'{e}te formula to write the identity (cf. \cite{jca})
$$\sum_{j=2}^{k} r_j (\lambda_j- \lambda_1) =0. \eqno(20)$$

{\bf Proposition 6}. {\it A possible non-trivial polynomial of degree $n \ge 2$ with $k$ distinct roots $\lambda_j,\ j=1,\dots, k$, sharing the  root $\lambda_1$ with its $n-1$st derivative must contain at least one root $\lambda_j,\ j\neq 1$ outside of the disk $D_\mu= \left\{ z \in \mathbb{C}:\  |z- \lambda_1-1| \le \mu \right\},\ \mu\in (0,1)$.} 

\begin{proof} The proof is based on the inequality involving weighted arithmetic and geometric means for complex numbers proved in \cite{four}.   Indeed, let $\lambda_j \in D_\mu,\ j\neq 1$.   Then according to \cite{four} we find from equality (20)  the following estimate 
$$0 = {1\over (1-\mu^2) (n-r_1)}\left| \sum_{j=2}^{k} r_j (\lambda_j- \lambda_1) \right| \ge 
{1-\mu\over 2\mu} \log\left( {1+\mu\over 1-\mu}\right) \exp \left( 1- {1-\mu\over 2\mu} \log\left( {1+\mu\over 1-\mu}\right) \right)$$

$$\times \prod_{j=2}^k  \left| \lambda_j- \lambda_1\right|^{r_j/(n-r_1)} \neq 0, $$
which  gives a contradiction. 

\end{proof}

Employing  the Sz.-Nagy type identity for complex roots  of $f$ and its $m$th derivative

$$   (z_{n-1} - z_{n-2})^2={1\over n(n-1)} \left[ \sum_{j=1}^k r_{j}(\lambda_j- z)^2-  n (z_{n-1} -  z)^2\right]$$
$$= {1\over (n-m)(n-m-1)} \left[ \sum_{j=1}^{n-m} (\xi_j^{(m)}- z)^2- ( n-m) (z_{n-1} -  z)^2\right],\  z \in \mathbb{C},\eqno(21)$$
which is proved in   \cite{jca}, we let $m=1$ and $z=0$, writing, in particular, the equality 
$$  \sum_{j=1}^{n-1} \left[\xi_j^{(1)}\right]^2= {n-2\over n} \sum_{j=1}^k r_{j} \lambda^2_j + z^2_{n-1}.\eqno(22)$$
The unique root $z_{n-1}$ of the $n-1$st derivative is called the centroid of the sets $\lambda_j$ and $\xi_j^{(m)}$  and it satisfies the  Sz.-Nagy type identity (cf. \cite{jca})
$$   z_{n-1} - z ={1\over n} \sum_{j=1}^k r_{j}(\lambda_j- z) = {1\over n-m }  \sum_{j=1}^{n-m} (\xi_j^{(m)}- z),\  z \in \mathbb{C}.\eqno(23)$$
But since among the roots  of the first derivative $f^\prime$ are  roots $\lambda_j$ of multiplicities $r_j-1$, correspondingly,  we let $m=1$ and $z=0$  in (23) to find
$$\sum_{j=1}^k  \lambda_j =  \sum_{j=1}^{k-1} \hat{\xi} _j^{(1)} +  z_{n-1}, \eqno(24)$$
where $\hat{\xi} _j^{(1)}$ are roots of the logarithmic derivative $\left(\log f(z) \right)^\prime$.  The corresponding identity for squares of these roots can be obtained from (21), (22),  and we have
$$  \sum_{j=1}^{k-1} \left[\hat{\xi}_j^{(1)}\right]^2=  \sum_{j=1}^k \lambda^2_j - 2(n-1) (z_{n-1} - z_{n-2})^2
- z_{n-1}^2.\eqno(25)$$
Following \cite{sch} we say that the set of $2n-1$ points is rectilinear, if $k$  roots $\lambda_j$ with multiplicities $r_j$ and $n-1$ roots $\xi_j^{(1)}$  are on a straight line in the complex plane, which passes through the origin.  Then the centroid $z_{n-1}$ and the root $z_{n-2}$ of the $n-2$nd derivative are contained on this line as well.  Hence there is an angle $\varphi$ such that  $\lambda_j = \pm |\lambda_j| e^{i\varphi},\    \xi_j^{(1)} = \pm |\xi_j^{(1)}| e^{i\varphi}$ for all $j$ and   $z_{n-1}=  \pm |z_{n-1}| e^{i\varphi},  z_{n-2}= \pm |z_{n-2}| e^{i\varphi} $.   Consequently,  equalities (22) and (25) imply the identities 
$$  \sum_{j=1}^{n-1} \left|\xi_j^{(1)}\right|^2= {n-2\over n} \sum_{j=1}^k r_{j} \left|\lambda\right| ^2_j + \left|z_{n-1}\right|^2,\eqno(26)$$
$$  \sum_{j=1}^{k-1} \left|\hat{\xi}_j^{(1)}\right|^2=  \sum_{j=1}^k \left|\lambda\right|^2_j - 2(n-1) \left( |z_{n-1}| - |z_{n-2}|\right)^2- \left|z_{n-1}\right|^2,\eqno(27)$$
respectively.  When  $z_{n-1}=  \pm |z_{n-1}| e^{i\varphi}$ but   $z_{n-2}= \mp  |z_{n-2}| e^{i\varphi} $ the latter equality becomes
$$  \sum_{j=1}^{k-1} \left|\hat{\xi}_j^{(1)}\right|^2=  \sum_{j=1}^k \left|\lambda\right|^2_j - 2(n-1) \left( |z_{n-1}| + |z_{n-2}|\right)^2- \left|z_{n-1}\right|^2.\eqno(28)$$
Furthermore, equality (26) suggests a generalization of the Schoenberg conjecture.  Indeed, we have 

{\bf Conjecture 1}. {\it For any complex roots of $f$ and its first derivative we have the inequality 
$$  \sum_{j=1}^{n-1} \left|\xi_j^{(1)}\right|^2\le  {n-2\over n} \sum_{j=1}^k r_{j} \left|\lambda\right|^2_j + \left|z_{n-1}\right|^2,\eqno(29)$$
with the equality sign if and only if all roots lie on a straight line, passing through the origin.} 

 The  Sz.-Nagy type identities (21) yield the following proposition for polynomials (1) whose centroid is zero.  
 
 {\bf Proposition 7.}  {\it  A monic polynomial $f$ of degree $n\ge 2$    whose roots lie on a straight line passing through the origin  is  $f(z)= z^n$ , if and only if  $z=0$  is the  double root of the $n-2$nd derivative.}

\begin{proof}  The necessity is obvious.  To prove the sufficiency we see that if $z_{n-2}= z_{n-1} =0$ then the first identity in  (21) with $z=0$ presumes 

$$ \sum_{j=1}^k r_{j}\lambda^2_j =0.$$
But since  $\lambda_j = \pm |\lambda_j| e^{i\varphi},\ j=1,\dots, k $ it has all roots should be zero and $f(z)= z^n$.
\end{proof} 

 Further, we generalize Proposition1 for polynomials whose roots are lying on the vertical or horizontal line of the complex plane.  In fact, we have

{\bf Proposition 8.}  {\it  A polynomial of degree $n\ge 2$    whose roots lie on the vertical (horizontal ) line of the complex plane  is trivial, if and only if its $n-2$nd derivative has a double root.}

\begin{proof}  Let distinct roots have the form $\lambda_j= a + i\ {\rm Im} \lambda_j,\  a \in \mathbb{R},\ j=1,\dots, k$ and $z_{n-2}$ is the double root of $n-2$nd derivative.  Then $z_{n-2}= z_{n-1}$. Hence the left-hand side of the first identity in (21) is zero.  Moreover, evidently, the centroid lies on the same vertical line. Letting  in (21)  $z= z_{n-1}$ it gives 
$$ \sum_{j=1}^k r_{j}(\lambda_j-  z_{n-1})^2 =0$$
or
$$\sum_{j=1}^k r_{j} {\rm Im} ^2\left[ \lambda_j-  z_{n-1}\right]  =0.$$ 
Hence,   ${\rm Im}  \lambda_j=  {\rm Im}  z_{n-1},\ j=2,\dots, k$ and the polynomial is trivial.  On the same manner we prove the proposition for roots lying on the horizontal line of the complex plane.  The necessity is obvious. 
\end{proof} 
Recalling identity (25) we generalize it for polynomials of degree $n\ge 2$ with simple $n$ roots $w_j,\  j= 1,\dots, n$ and roots of its $m$-th derivative.  Precisely, with the use of (21) we find
$$  \sum_{j=1}^{n-m} \left[\xi_j^{(m)}\right]^2 =  \sum_{j=1}^n w_j^2  - m(2n-m-1) (z_{n-1} - z_{n-2})^2
- m  z^2_{n-1},\  m=0,1,\dots, n.\eqno(30)$$
Next    we will derive a formula, involving higher order derivatives of $\log f(z)$, which seems to be new. 

{\bf Lemma 2}. {\it Let $m \in \mathbb{N}_0,\ z \in \mathbb{C}$ and $f(z)\neq 0$. Then the following formula takes place}

$$\left( \log f(z)\right)^{(m+1)} =  \left( {f^\prime(z)\over f(z)} \right)^{(m)}  
= \sum_{j=0}^m {(-1)^j\over j+1} \   {{m+1}\choose {j+1}}  \frac {\left( \left[ f(z)\right] ^{j+1} \right)^{(m+1)} }{  \left[ f(z)\right] ^{j+1} }.\eqno(31)$$

\begin{proof}  In order to prove (31) we call the familiar Hoppe formula  (see \cite{jon}, p. 224) for higher derivatives of the composition of two functions.  Thus we derive
$$\left( \log f(z)\right)^{(m+1)} =  \left( {f^\prime(z)\over f(z)} \right)^{(m)}  
= \sum_{s=0}^{m+1}  {(-1)^s (s-1)! \over s!  \left[ f(z)\right] ^{s} }  \sum_{j =0}^{s}  (-1)^{s-j}  {{s}\choose {j}}  
\left[ f(z)\right] ^{s-j} \left( \left[ f(z)\right] ^{j} \right)^{(m+1)} $$
$$ = \sum_{s= 0}^{m+1}  (s-1) !    \sum_{j =0}^{s}  { (-1)^{j+1} \over j! (s -j)! }   
\frac{ \left( \left[ f(z)\right] ^{j} \right)^{(m+1)}}{ \left[ f(z)\right] ^{j} } $$
$$ =    \sum_{j =0}^{m+1}   (-1)^{j+1}  \frac{ \left( \left[ f(z)\right] ^{j} \right)^{(m+1)}}{ j!\  \left[ f(z)\right] ^{j} } \sum_{s= j-1}^{m}  {s! \over  (s+1 -j)! } $$
$$ =    \sum_{j =1}^{m+1}   (-1)^{j+1}  \frac{ \left( \left[ f(z)\right] ^{j} \right)^{(m+1)}}{ j!\  \left[ f(z)\right] ^{j} } \sum_{s= j-1}^{m}  {s! \over  (s+1 -j)! }. $$
Hence a simple substitution in the index of summation  and the use of the combinatorial  identity 
$$\sum_{s= j}^{m} {{s}\choose {j}}   =  {{m+1}\choose {j+1}}$$
lead to (31) and complete the proof of Lemma 2.    

\end{proof}

Let a possible non-trivial CA- polynomial $f$ have real zeros only .  Proposition $5$ says that  it has at least $5$ distinct zeros, i.e $k\ge 5$. Then, by virtue of the Rolle theorem all zeros of the derivatives $f^{(j)} (x),\ x \in \mathbb{R}, \  j =r-1, r,\dots, n-1$, where  $2 \le r= \max_{1\le j\le k} r_j$  are simple.  Denoting by 
$$A= \{ \lambda_1, \dots, \lambda_k \},\quad  B_j =  \{  \xi_1^{(j)},\dots, \xi_{n-j}^{(j)} \}, \ j=r-1, r,\dots, n-1$$
sets of roots of $f$ and its $m$th derivatives,  we have  by definition of the CA-polynomial $ B_j \cap A \neq \emptyset,\  j =r-1, r,\dots, n-1$.   Moreover, letting in (31) $m=1$,   we write the Laguerre inequality for derivatives $f^{(j)}$ (see, for instance, in \cite{bor}) 
$$ {d^2\over dx^2} \left[\log f^{(j)} (x) \right]  < 0,\ j= r-1, r,\dots, n-2,\eqno(32)$$
or, in the equivalent  form,  
$$  \left[ f^{(j+1)} (x) \right]^2 >   f^{(j)} (x) f^{(j+2)} (x),\  j= r-1, r,\dots, n-2.\eqno(33)$$
The latter inequality implies the property  $B_j \cap B_{j+1}  = \emptyset.$  Let  $C_j \subset A$ be a  subset of $B_j$, containing $n_j\in \mathbb{N}$ common roots of $f^{(j)}$ with $f$, i.e.
$$C_j= \{  \lambda_{j,1},\dots, \lambda_{j, n_j} \} \subset  B_j.$$
Clearly, the number $n_{r-1} $ of common roots with the $r-1$st derivative does not exceed $\min (k-1,\  n-r+1)$ and $n_j \le \min (k-2,\ n-j) \ j= r,\dots, n-1$ because the minimal and maximal roots of $f$ cannot be zeros of $f^{(j)},\  j \ge r.$ 
Condition (33) says that it may happen that $C_j \cap C_{j+2}  \neq  \emptyset,\  2\le s \le n- r.$ Writing equality (30) for roots of the $m$ th and $m+s$ th derivatives, where $r-1 \le m \le n-s-1$, we find 
$$  \sum_{j=1}^{n-m} \left[\xi_j^{(m)}\right]^2-    \sum_{j=1}^{n-m-s} \left[\xi_j^{(m+s)}\right]^2=   s(2(n-m) -s-1) (x_{n-1} - x_{n-2})^2 + s x^2_{n-1}.\eqno(34)$$
Hence we arrive at

{\bf Proposition 9.}  {\it  Let $m= r-1, r,\dots, n-1$, where  $2 \le r= \max_{1\le j\le k} r_j $ and $2\le s \le n- r.$  Then roots  of the $m$ th and $m+s$th derivatives of a  possible non-trivial CA-polynomial  with only real zeros satisfy the condition }
$$  \sum_{1\le j\le n-m, \  \xi_j^{(m)} \notin \ C_m \cap C_{m+s} } \left[\xi_j^{(m)}\right]^2 \ge   
\sum_{1\le j\le n-m-s, \  \xi_j^{(m)} \notin \  C_m \cap C_{m+s}} \left[\xi_j^{(m+s)}\right]^2.$$

\bibliographystyle{amsplain}

\end{document}